\documentclass{article}

\usepackage{graphicx}

\newtheorem{thm}{Theorem}[section]
\newtheorem{cor}[thm]{Corollary}

\newtheorem{de}[thm]{Definition}

\begin{document}

\title{The Shape Parameter in the Gaussian Function}         
\author{Lin-Tian Luh\\Department of Mathematics, Providence University\\Shalu, Taichung, Taiwan\\
Email:ltluh@pu.edu.tw\\Fax:886-4-26324653, Tel:886-4-26328001ext.15126 }        
\date{\today}          
\maketitle
{\bf Abstract}. This is the fifth of our series of works about the shape parameter. We now explore the parameter $\beta$ contained in the famous gaussian function $e^{-\beta |x|^{2}},\ x\in R^{n}$. In the theory of radial basis functions(RBF), gaussian is frequently used in virtue of its good error bound and numerical tractability. However the optimal choice of $\beta$ has been unknown. RBF people only know that $\beta$ is very influential, but do not have a reliable criterion of its choice. The purpose of this paper is to uncover its mystery.
\section{Introduction}  
Let
\begin{equation}
  h(x):=e^{-\beta |x|^{2}},\ x\in R^{n},\ \beta >0.
\end{equation}
For any given data points $(x_{j},y_{j}),\ j=1,\cdots,N$, where $X=\{x_{1},\cdots,x_{N}\}$ is a subset of $R^{n}$ and the $y_{j}'s$ are real or complex numbers, it's well known that there is a so-called $h$ spline interpolant $s(x)$ of these data points, defined by
\begin{equation}
  s(x)=p(x)+\sum_{j=1}^{N}c_{j}h(x-x_{j})
\end{equation}    
where $p(x)$ is a polynomial in $P_{m-1}$, $m$ being the order of the conditional positive definiteness of $h$, and $c_{j}'s$ are chosen so that
\begin{equation}
  \sum_{j=1}^{N}c_{j}q(x_{j})=0
\end{equation}
for all polynomials $q$ in $P_{m-1}$ and 
\begin{equation}
  p(x_{i})+\sum_{j=1}^{N}c_{j}h(x_{i}-x_{j})=y_{i},\ i=1,\cdots,N.
\end{equation}
Here $P_{m-1}$ denotes the class of polynomials in $R^{n}$ of degree $\leq m-1$, and $P_{m-1}:=\{0\}$ if $m=0$. The linear system of equations (3) and (4) has a unique solution $(c_{1},\cdots,c_{N})$ with $p(x)\equiv 0$ since $m=0$ for (1). Further details can be found in \cite{MN1}.

We use $s(x)$ as the approximating function. For error estimates there are two kinds of error bound, the algebraic type \cite{MN2} and the exponential type \cite{MN3} and \cite{Lu3-1}. The latter is an improved form of the former. In this paper we will show that if the shape parameter $\beta$ is chosen in an appropriate way, the exponential-type error bound will be minimal.
\section{Fundamental Theory}
In the theory of RBF, any conditionally positive definite(c.p.d.) function $h$ of order $m$ induces a function space called {\bf native space}, denoted by {\bf ${\cal C}_{h,m}$}. Its definition and characterization can be found in \cite{Lu1}, \cite{Lu2}, \cite{Lu4}, \cite{MN1}, \cite{MN2} and \cite{We}. We adopt Madych and Nelson's definition which is quite different from that of Wu and Schaback\cite{WS}. Besides the native space, each $f\in {\cal C}_{h,m}$ has a seminorm $\|f\|_{h}$ which plays an important role in our theory. In this paper $h$ is just the function in (1) and $m=0$.

There are two main theorems in this paper. Before introducing the first of our main theorems, we need a basic definition.    
\begin{de}
  For any positive integer n, the number $\gamma_{n}$ is defined by $\gamma_{1}=2$ and $\gamma_{n}=2n(1+\gamma_{n-1})$ if $n\geq 2$. 
\end{de}
As pointed out in section1, the theoretical cornerstone of our approach is the exponential-type error bound constructed by Madych and Nelson. However some crucial constants in the error bound had been unknown and considered to be a hard question. Fortunately these constants are thoroughly clarified in \cite{Lu3}. In \cite{Lu3} the author presents a complete and lucid exponential-type error bound for gaussian interpolation. In order to develop useful criteria of the optimal choice of $\beta$, we need the following theorem which is cited from \cite{Lu3} directly.
\begin{thm}
Let $h(x)=e^{-\beta \left\vert x\right\vert ^{2}},$ $\beta >0,$ be the
gaussian function in $R^{n}.$ Then, given a positive number $b_{0},$ there are positive
constants $\delta _{0},c,$and $C$ for which the following is true: If $f\in 
{\cal C}_{h,m}$ and $s$ is the $h$ spline that interpolates $f$ on a
subset $X$ of $R^{n},$ then
\begin{equation}
\left\vert f(x)-s(x)\right\vert \leq \triangle ^{^{\prime \prime
}}(C_{\delta })^{\frac{c}{\delta }}\cdot \left\Vert f\right\Vert _{h}
\end{equation}
, where $\Delta ''=\pi^{\frac{n-1}{4}}\cdot (n\cdot \alpha_{n})^{\frac{1}{2}}\cdot 2^{\frac{n+1}{4}}\cdot \left( \frac{\sqrt{3}}{e}\right) ^{\frac{n-2}{4}}$ for even n and $\Delta ''=\pi^{\frac{n}{4}}\cdot (n\cdot \alpha_{n})^{\frac{1}{2}}\cdot \left( \frac{\sqrt{3}}{e}\right) ^{\frac{n-1}{4}}$ for odd n, holds for all x in a cube E provided that (a)E has side b and b$\geq b_{0}$,(b)$0<\delta \leq \delta_{0}$, and (c)every subcube of E of side $\delta$ contains a point of X. Here, $\alpha_{n}$ denotes the volume of the unit ball in $R^{n}$.

The number c is equal to $\frac{b_{0}}{8\gamma_{n}}$ where $\gamma_{n}$ was defined in Definition2.1. The number C is equal to $(3^{\frac{3}{4}}\cdot e\cdot \sqrt{2\rho \beta}\sqrt{n}\cdot e^{2n\gamma_{n}})^{4}\cdot b_{0}^{3}\cdot \gamma_{n}$, where $\rho=\frac{\sqrt{3}}{e}$. Moreover, $\delta_{0}$ can be defined by
$$\delta_{0}=\min \left\{ \frac{1}{(3^{\frac{3}{4}}\cdot e\cdot \sqrt{2\rho \beta}\cdot \sqrt{n}\cdot e^{2n\gamma_{n}})^{4}\cdot b_{0}^{3}\cdot \gamma_{n}},\ \delta_{n}\ \right\}$$
, where
$$\delta_{n}=\left\{ \begin{array}{llll}
                      \frac{b_{0}}{2\gamma_{n}},      & \mbox{if $n=1$,}\\
                      \frac{b_{0}}{2\gamma_{n}(n-1)}, & \mbox{if n is odd and $n>1$,}\\
                      \frac{b_{0}}{2\gamma_{n}(n-2)},& \mbox{if n is even and $n>2$,}\\
                      \frac{b_{0}}{2\gamma_{n}},      & \mbox{if $n=2$.}
                    \end{array}  \right.$$
\end{thm}
Note that the error bound (5) tends to zero as $\delta$ tends to zero. This is the key to understanding this seemingly complicated theorem. In (5) the constant $C$ highly depends on $\beta$. It's tempting to think that in (5) only $C$ is influenced by $\beta$. In fact, $\|f\|_{h}$ also changes as $\beta$ changes. The change of $\|f\|_{h}$ cannot be seen in a transparent way. Therefore Theorem2.2 cnnot be used directly to find the optimal $\beta$.

In order to overcome this problem, we introduce two function spaces as follows.
\begin{de}
  For any $\sigma>0$, the class of band-limited functions in $L^{2}(R^{n})$ is
$$B_{\sigma}:=\{f\in L^{2}(R^{n}):\hat{f}(\xi)=0\ if\ |\xi|>\sigma \}$$
, where $\hat{f}$ denotes the Fourier transform of $f$.
\end{de}
\begin{de}
For any $\sigma>0$,
$$G_{\sigma}:=\{f\in L^{2}(R^{n}):\ \int|\hat{f}(\xi)|^{2}e^{\frac{|\xi|^{3}}{\sigma}}d\xi <\infty \}$$
, where $\hat{f}$ denotes the Fourier transform of $f$. For each $f\in G_{\sigma}$,
$$\|f\|_{G_{\sigma}}:=\left\{ \int |\hat{f}(\xi)|^{2}e^{\frac{|\xi|^{3}}{\sigma}}d\xi \right\} ^{1/2}.$$
\end{de}
Let's investigate $B_{\sigma}$ first. For each $f\in B_{\sigma}$, by Corollary3.3 of \cite{MN2},
\begin{eqnarray*}
  \|f\|_{h} & = & \left\{ \int\frac{|\hat{f}(\xi)|^{2}}{(2\pi)^{2n}\hat{h}(\xi)}d\xi\right\}^{1/2}\\
        & = & \frac{1}{(2\pi)^{n}}\left\{\int|\hat{f}(\xi)|^{2}e^{\frac{\|\xi\|^{2}}{4\beta}}d\xi\right\}^{1/2}\\
        & \leq & \frac{1}{(2\pi)^{n}}e^{\frac{\sigma^{2}}{8\beta}}\left\{\int|\hat{f}(\xi)|^{2}d\xi\right\}^{1/2}\\
        & = & \frac{1}{(2\pi)^{n}}e^{\frac{\sigma^{2}}{8\beta}}\|f\|_{L^{2}(R^{n})}. 
\end{eqnarray*}
Substituting this result into (5), we obtain the following useful theorem.
\begin{thm}
  For any $\sigma>0$, $f\in B_{\sigma}$ implies $f\in {\cal C}_{h,m}$ and (5) can be transformed into 
\begin{equation}
  |f(x)-s(x)|\leq \Delta''(C\delta)^{\frac{c}{\delta}}\cdot (2\pi)^{-n}\cdot e^{\frac{\sigma^{2}}{8\beta}}\cdot \|f\|_{L^{2}(R^{n})}
\end{equation}
, where $h$ is defined as in (1).
\end{thm}
Functions in $G_{\sigma}$ can be treated in a similar way. For any $f\in G_{\sigma}$, we have
\begin{eqnarray*}
  \|f\|_{h} & = & \frac{1}{(2\pi)^{n}}\left\{\int|\hat{f}(\xi)|^{2}e^{\frac{|\xi|^{2}}{4\beta}}d\xi\right\}^{1/2}\\
            & = & \frac{1}{(2\pi)^{n}}\left\{\int|\hat{f}(\xi)|^{2}e^{\frac{|\xi|^{3}}{\sigma}}e^{\frac{|\xi|^{2}}{4\beta}-\frac{|\xi|^{3}}{\sigma}}d\xi\right\}^{1/2}\\
            & \leq & \frac{1}{(2\pi)^{n}}\left\{ \sup_{\xi\in R^{n}}e^{\frac{|\xi|^{2}}{4\beta}-\frac{|\xi|^{3}}{\sigma}}\right\}^{1/2}\|f\|_{G_{\sigma}}.
\end{eqnarray*}
This gives the following theorem.
\begin{thm}
  For any $\sigma>0,\ f\in G_{\sigma}$ implies $f\in {\cal C}_{h,m}$ and (5) can be transformed into
\begin{equation}
  |f(x)-s(x)|\leq \Delta''(C\delta)^{\frac{c}{\delta}}(2\pi)^{-n}\left\{ \sup_{\xi\in R^{n}}e^{\frac{|\xi|^{2}}{4\beta}-\frac{|\xi|^{3}}{\sigma}}\right\}^{1/2}\|f\|_{G_{\sigma}}
\end{equation}
, where $h$ is defined as in (1).
\end{thm}
There is an improved exponential-type error bound from which a set of criteria for the choice of the shape parameter can also be developed. In this kind of error bound the data points are not purely scattered and the interpolation happens in an n-simplex.

The definition of n-simplex in $R^{n}$ can be found in \cite{Fl}. 1-simplex is a line segment. 2-simplex is a triangle. 3-simplex is a tetrahedron.

Let $T_{n}$ be an n-simplex in $R^{n}$ with vertices $v_{1},\cdots,v_{n+1}$. Any $x\in T_{n}$ can be written as a convex combination of the vertices:
$$x=\sum_{i=1}^{n+1}c_{i}v_{i}$$
where $\sum_{i=1}^{n+1}c_{i}=1$ and $c_{i}\geq 0$ for all $i$. We call $(c_{1},\cdots,c_{n+1})$ the barycentric coordinate of $x$. Let's define {\bf `evenly spaced' points} of degree $l$ to be those points whose barycentric coordinates are of the form
$$(\frac{k_{1}}{l},\frac{k_{2}}{l},\cdots,\frac{k_{n+1}}{l}),\ k_{i}\ nonnegative\ integers\ with\ \sum_{i=1}^{n+1}k_{i}=l.$$

It's easily seen that the number of such points in $T_{n}$ is exactly $N=dimP_{l}^{n}$. Also, as stated in \cite{Bo}, evenly spaced points form a determining set for $P_{l}^{n}$, the space of polynomials of degree $l$ in n variables.

We can base our criteria of choosing $\beta$ on a crucial theorem which we cite directly from \cite{Lu3-1} with a slight modification.
\begin{thm}
  Let $h(x):=e^{-\beta|x|^{2}}$ be the gaussian function in $R^{n}$. For any positive number $b_{0}$, there are positive constants $\delta_{0},c_{1},c_{2}$, and $c_{3}$ independent of $n$, for which the following is true: If $f\in {\cal C}_{h,m}$ , the native space induced by $h$, and $s$ is the $h$ spline that interpolates $f$ on a subset $X$ of $R^{n}$, then 
\begin{equation}
  |f(x)-s(x)|\leq c_{1}\sqrt{\delta}(c_{2}\delta)^{\frac{c_{3}}{\delta}}\cdot \| f\| _{h}
\end{equation}
for all $x$ in a subset $\Omega$ of $R^{n}$, and $0<\delta \leq \delta_{0}$, where $\Omega$ satisfies the property that for any $x$ in $\Omega$ and any number $\frac{b_{0}}{2}\leq r\leq b_{0}$, there is an $n$ simplex $Q$ with diameter $diamQ=r,x\in Q\subseteq \Omega$, such that for any integer $l$ with $\frac{b_{0}}{\delta}\leq l\leq \frac{2b_{0}}{\delta}$, there is on $Q$ an evenly spaced set of centers from $X$ of degree $l-1$.(In fact, the set $X$ can be chosen to consist of these evenly spaced centers in $Q$ only.) Here $\|f\| _{h}$ is the $h$-norm of $f$ in the native space. The numbers $\delta_{0},c_{1},c_{2}$, and $c_{3}$ are given by

$ \delta_{0}:=\min \left\{ b_{0},\ \frac{1}{\rho_{3}^{4}\cdot 3^{3}\cdot 2^{7}\cdot b_{0}^{3}}\right\} \ \ where\ \ \rho_{3}=12^{\frac{1}{4}}\cdot \sqrt{e\beta}\ \ $; 
$\left\{
\begin{array}{lll}
 c_{1}:=\left\{ \begin{array}{ll}
                                                                                          \Delta''\cdot \frac{1}{\sqrt{16\pi}}\cdot \frac{1}{\sqrt{b_{0}}}\ \ for\ \ odd\ \ n, \\ \Delta''\cdot \frac{1}{\sqrt{16\pi}}\cdot \frac{1}{\sqrt{b_{0}}}\ \ for\ \ even\ \ n,
                \end{array}\right. \  \\ c_{2}:=\rho_{3}^{4}\cdot 3^{3}\cdot 2^{7}\cdot b_{0}^{3}\ \ , \\ c_{3}:=\frac{b_{0}}{4},
 \end{array} \right.$
 where $\Delta''$ is defined by

$\Delta'':=\left\{ \begin{array}{ll}
                     \sqrt{2+\frac{1}{e}}\cdot \pi^{\frac{n-1}{4}}\cdot (n\alpha_{n})^{\frac{1}{2}}\cdot 2^{\frac{n}{4}}\cdot \rho^{\frac{n-1}{4}}\ for\ odd\ n, \\ \pi^{\frac{n-1}{4}}\cdot (n\alpha_{n})^{\frac{1}{2}}\cdot 2^{\frac{n+1}{4}}\cdot \rho^{\frac{n-2}{4}}\ for\ even\ n,\ with\ \rho=\frac{\sqrt{3}}{e}   
                   \end{array} \right. $ ,where the number $\alpha_{n}$ denotes the volume of the unit ball in $R^{n}$. 
In particular, if the point $x$ in $\Omega$ is fixed, the only requirement for $\Omega$ is the existence of an n simplex $Q$, with $diamQ=r,\ x\in Q\subseteq \Omega$, satisfying the afore-mentioned property of evenly spaced centers, and the centers $x_{1},\cdots, x_{N}$ in $s(x)$, as defined in (2), are just the evenly spaced points in $Q$.
\end{thm}
{\bf Remark}: This seemingly complicated theorem is in fact not difficult to understand. The number $\delta$ is in spirit equivalent to the well known fill-distance. The error bound (8) tends to zero rapidly as $\delta$ tends to zero. In practical application only finitely many interpolation points will be involved. Hence only a finite number of simplices will appear, even if $\Omega$ is unbounded. Note that among $c_{1},c_{2}$ and $c_{3}$, only $c_{2}$ depends on the shape parameter $\beta$. It's tempting to think that the influence of $\beta$ on (8) happens only in $c_{2}$. This is wrong. Be careful that, as in (5), $\|f\|_{h}$ depends on $\beta$ also.

By transforming the norms we obtain the following useful results now.

\begin{cor}
  Suppose $f\in B_{\sigma},\sigma>0$. Then (8) can be transformed into 
\begin{equation}
  |f(x)-s(x)|\leq c_{1}\sqrt{\delta}(c_{2}\delta)^{\frac{c_{3}}{\delta}}\cdot (2\pi)^{-n}\cdot e^{\frac{\sigma^{2}}{8\beta}}\cdot \|f\|_{L^{2}(R^{n})}.
\end{equation}
\end{cor}

\begin{cor}
  Suppose $f\in G_{\sigma},\sigma>0$. Then (8) can be transformed into
\begin{equation}
  |f(x)-s(x)|\leq c_{1}\sqrt{\delta}(c_{2}\delta)^{\frac{c_{3}}{\delta}}\frac{1}{(2\pi)^{n}}\left\{\sup_{\xi\in R^{n}}e^{\frac{|\xi|^{2}}{4\beta}-\frac{|\xi|^{3}}{\sigma}}\right\}^{1/2}\|f\|_{G_{\sigma}}
\end{equation}.
\end{cor}
\section{Criteria of Choosing $\beta$--the Scattered Type}
In (6) the parts influenced by $\beta$ are $(C\delta)^{\frac{c}{\delta}}$ and $e^{\frac{\sigma^{2}}{8\beta}}$. By the definition of $C$ and $c$, one can easily find that the function
$$MN(\beta):=\beta^{\frac{b_{0}}{4\gamma_{n}\delta}}\cdot e^{\frac{\sigma^{2}}{8\beta}}$$
describes the dependence of the error bound (6) on $\beta$. Let's call this function the MN function and its graph the MN curve. Obviously, finding the optimal $\beta$ is equivalent to finding the value $\beta$ minimizing $MN(\beta)$.

Similarly. in (7) the MN function is
$$MN(\beta):=\beta^{\frac{b_{0}}{4\gamma_{n}\delta}}\left\{sup_{\xi\in R^{n}}e^{\frac{|\xi|^{2}}{4\beta}-\frac{|\xi|^{3}}{\sigma}}\right\}^{1/2}.$$

Another important thing which remains to be investigated is $\delta$. Our fundamental theory begins at Theorem2.2 where $0<\delta\leq \delta_{0}$ is a requirement. By the definition of $\delta_{0}$, we find
\begin{eqnarray*}
  \delta\leq \delta_{0} & iff & \delta\leq \delta_{n}\ and\ \delta\leq \frac{1}{(3^{\frac{3}{4}}e\sqrt{2\rho\beta n}e^{2n\gamma_{n}})^{4}b_{0}^{3}\gamma_{n}}\\
                        & iff & \delta\leq \delta_{n}\ and\ \beta^{2}\leq \frac{1}{(3^{\frac{3}{4}}e\sqrt{2\rho n}e^{2n\gamma_{n}})^{4}b_{0}^{3}\gamma_{n}\delta}\\
                        & iff & \delta\leq \delta_{n}\ and\ \beta\leq \frac{1}{3^{\frac{3}{2}}2\rho n e^{4n\gamma_{n}+2}b_{0}^{\frac{3}{2}}\sqrt{\gamma_{n}}\sqrt{\delta}}.  
\end{eqnarray*}
This severely restricts the ranges of $\delta$ and $\beta$. However the upper bound of $\beta$ can be made arbitrarily large because $\delta$ can be arbitrarily small, theoretically.

We summarize these results in the following two criteria.\\
\\ 
{\bf Case1}. \fbox{$f\in B_{\sigma}$} Let $\sigma>0$ and $f\in B_{\sigma}$. Under the conditions of Theorem2.2, for any fixed $b_{0}>0$ and $\delta\leq \delta_{0}$, the optimal value of $\beta$ in the interval $(0,\beta_{0}]$ where $\beta_{0}= \frac{1}{3^{\frac{3}{2}}2\rho n e^{4n\gamma_{n}+2}b_{0}^{\frac{3}{2}}\sqrt{\gamma_{n}}\sqrt{\delta}}$ is the number minimizing
$$MN(\beta):=\beta^{\frac{b_{0}}{4\gamma_{n}\delta}}e^{\frac{\sigma^{2}}{8\beta}}.$$
{\bf Case2}. \fbox{$f\in G_{\sigma}$} Let $\sigma>0$ and $f\in G_{\sigma}$. Under the conditions of Theorem2.2, for any fixed $b_{0}>0$ and $\delta\leq \delta_{0}$, the optimal value of $\beta$ in the interval $(0, \beta_{0}]$ where $\beta_{0}=\frac{1}{3^{\frac{3}{2}}2\rho n e^{4n\gamma_{n}+2}b_{0}^{\frac{3}{2}}\sqrt{\gamma_{n}}\sqrt{\delta}}$ is the number minimizing
$$MN(\beta):=\beta^{\frac{b_{0}}{4\gamma_{n}\delta}}\left\{\sup_{\xi\in R^{n}}e^{\frac{|\xi|^{2}}{4\beta}-\frac{|\xi|^{3}}{\sigma}}\right\}^{1/2}$$\\
\\
{\bf Remark}:(a)In both cases $MN(\beta)\rightarrow \infty$ both as $\beta\rightarrow 0^{+}$ and $\beta\rightarrow \infty$. (b)The number $\beta$ minimizing $MN(\beta)$ can be obtained by Mathematica or Matlab.\\
\\
{\bf Numerical Results}:We now present some pictures of the MN function. The lowest point of the MN curve corresponds to the optimal choice of the shape parameter $\beta$. In order to make the graph look better, we sometimes multiply the function by a constant k and call it the modified MN function. Its graph will then be called the modified MN curve.

\begin{figure}[h]
\centering
\includegraphics[scale=1.0]{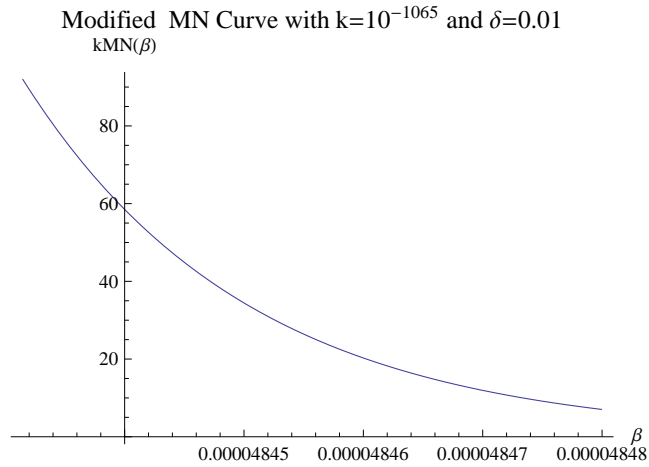}
\caption{$f\in B_{\sigma},n=1,\sigma=1,b_{0}=1$}
\end{figure}

\clearpage

\begin{figure}[t]
\centering
\includegraphics[scale=1.0]{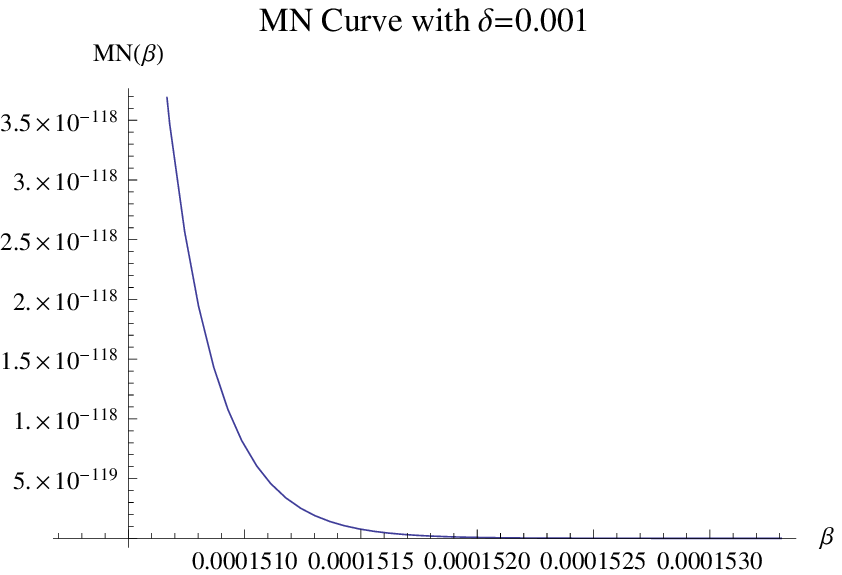}
\caption{$f\in B_{\sigma},n=1,\sigma=1,b_{0}=1$}

\includegraphics[scale=1.0]{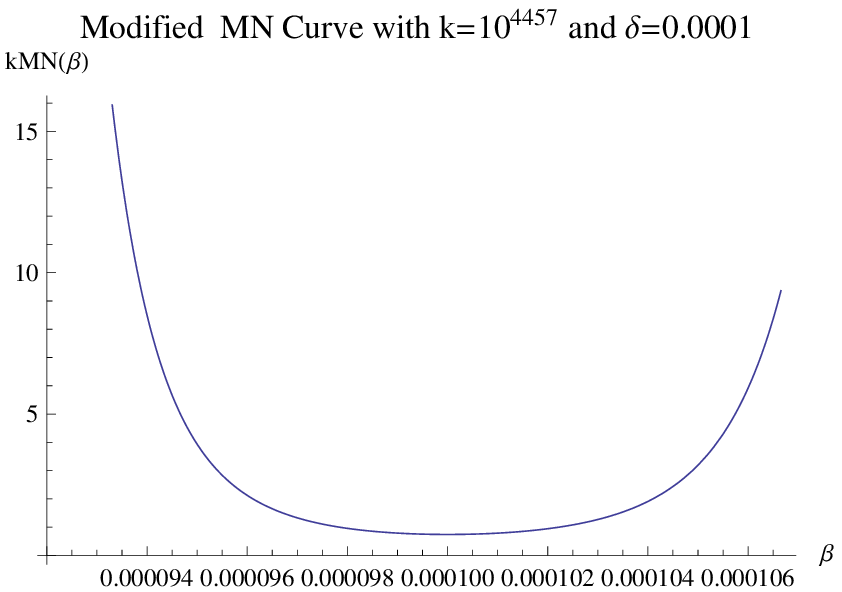}
\caption{$f\in B_{\sigma},n=1,\sigma=1,b_{0}=1$}

\includegraphics[scale=1.0]{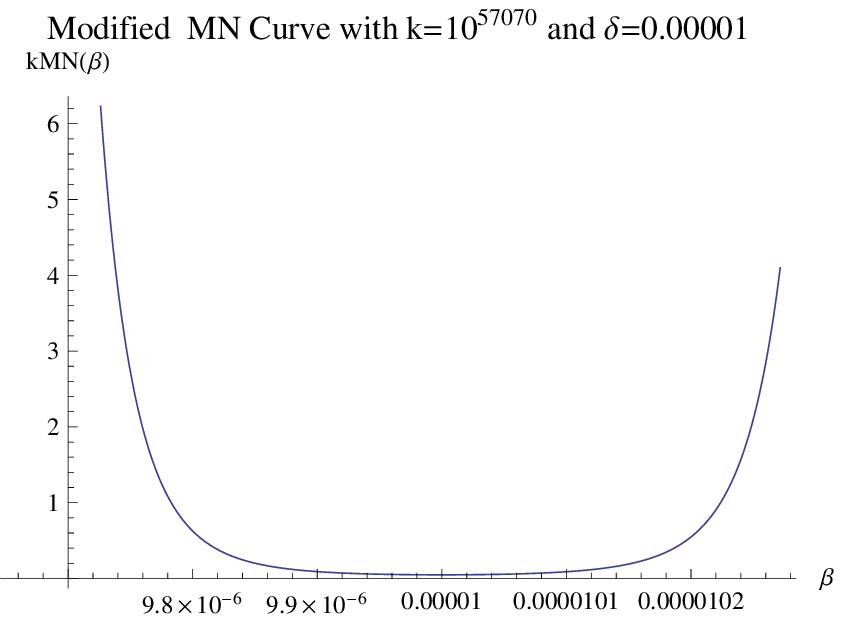}
\caption{$f\in B_{\sigma},n=1,\sigma=1,b_{0}=1$}

\end{figure}
\clearpage

\begin{figure}[t]
\centering
\includegraphics[scale=1.0]{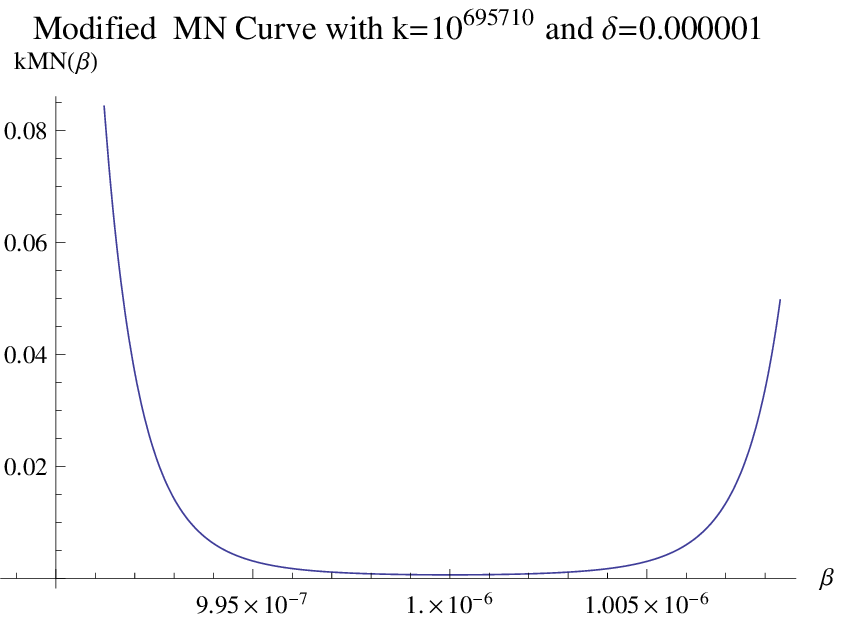}
\caption{$f\in B_{\sigma},n=1,\sigma=1,b_{0}=1$}

\includegraphics[scale=1.0]{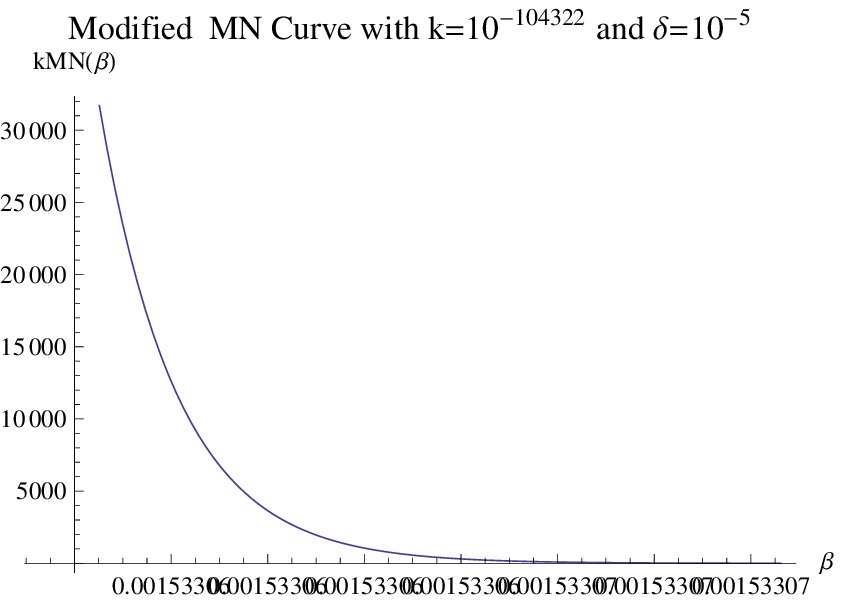}
\caption{$f\in G_{\sigma},n=1,\sigma=1,b_{0}=1$}

\includegraphics[scale=1.0]{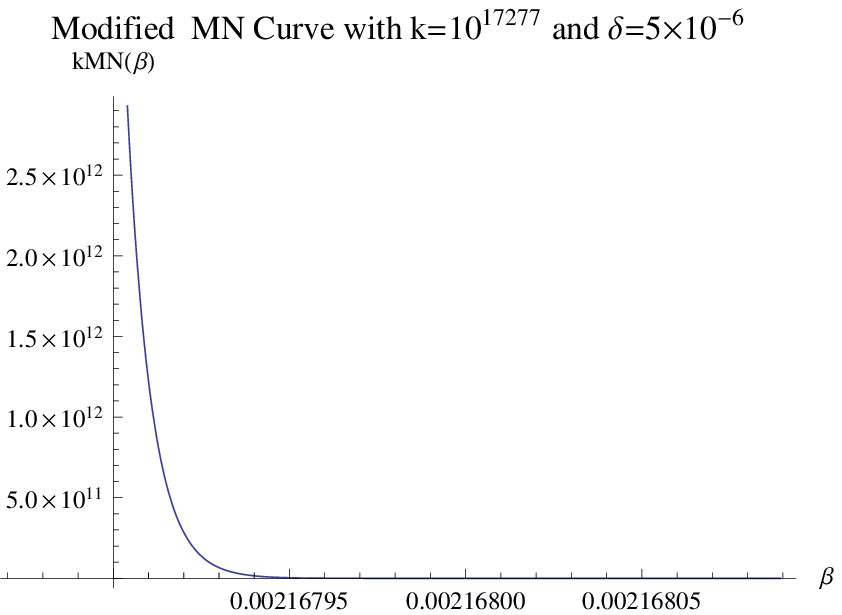}
\caption{$f\in G_{\sigma},n=1,\sigma=1,b_{0}=1$}

\end{figure}

\clearpage

\begin{figure}[t]
\centering
\includegraphics[scale=1.0]{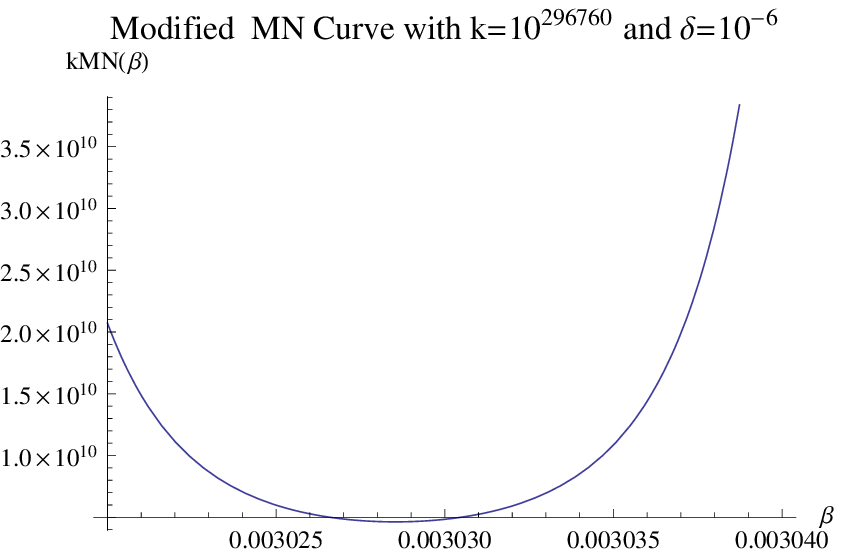}
\caption{$f\in G_{\sigma},n=1,\sigma=1,b_{0}=1$}

\includegraphics[scale=1.0]{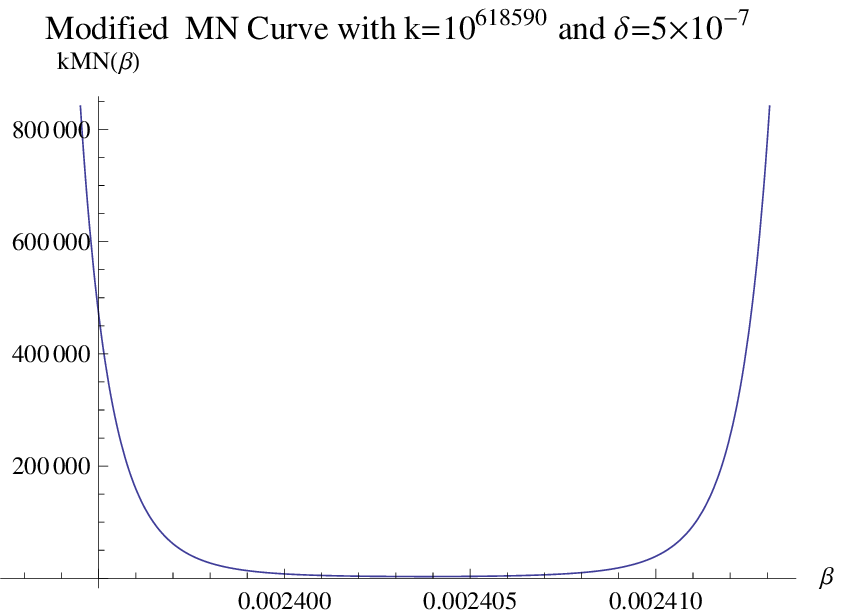}
\caption{$f\in G_{\sigma},n=1,\sigma=1,b_{0}=1$}

\includegraphics[scale=1.0]{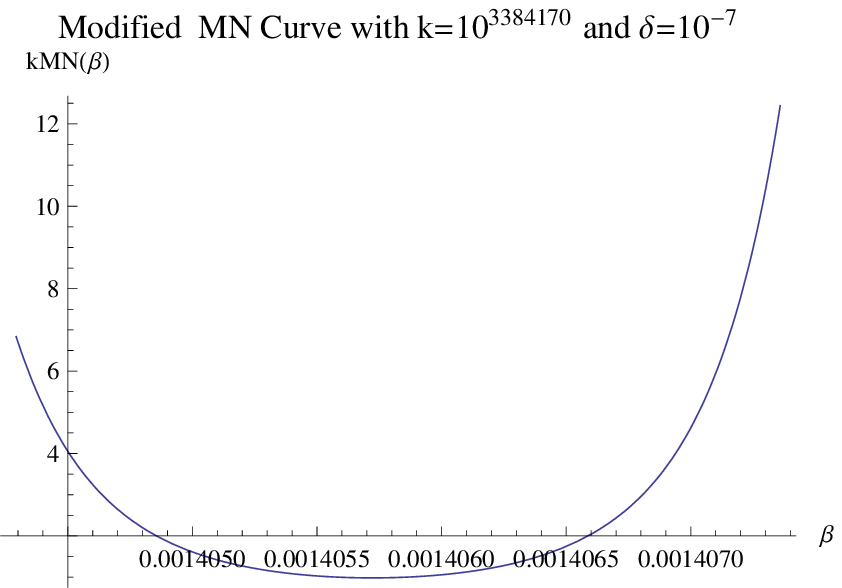}
\caption{$f\in G_{\sigma},n=1,\sigma=1,b_{0}=1$}
\end{figure}

\clearpage

Although we didn't show the whole curves of the MN or modified MN functions, the crucial parts were presented.\\
\\
{\bf Remark}: In our criteria both $\delta$ and $\beta$ have an upper bound which is quite restrictive, especially for high dimensions. This seems to be an important topic and deserves future research.

\section{Criteria of Choosing $\beta$--the Evenly Spaced Type}
Formula (9) and (10) provide us with a very good theoretical ground to choose $\beta$.

In (9), if we extract the parts influenced by $\beta$, we will get a function of $\beta$, i.e.
$$(\beta^{2})^{\frac{b_{0}}{4\delta}}e^{\frac{\sigma^{2}}{8\beta}}.$$
As before, let's call it the MN function and denote it by $MN(\beta)$. Thus
$$MN(\beta):=\beta^{\frac{b_{0}}{2\delta}}e^{\frac{\sigma^{2}}{8\beta}}.$$
Its graph will be called the MN curve. Finding the optimal $\beta$ is then equivalent to finding the number minimizing $MN(\beta)$. 
Similarly, in (10), the MN function is
$$MN(\beta):=\beta^{\frac{b_{0}}{2\delta}}\left\{\sup_{\xi\in R^{n}}e^{\frac{|\xi|^{2}}{4\beta}-\frac{|\xi|^{3}}{\sigma}}\right\}^{1/2}.$$
Note that in Theorem2.7 we require that $\delta\leq \delta_{0}$ where $\delta_{0}$ highly depends on $\beta$ and should be treated in a rigorous way.

By the definition of $\delta_{0}$,
\begin{eqnarray*}
  \delta\leq \delta_{0} & iff & \delta\leq b_{0}\ and\ \delta\leq \frac{1}{12e^{2}\beta^{2}27\cdot2^{7}b_{0}^{3}}\\
                        & iff & \delta\leq b_{0}\ and\ \beta^{2}\leq \frac{1}{12e^{2}27\cdot 2^{7}b_{0}^{3}\delta}\\
                        & iff & \delta\leq b_{0}\ and\ \beta\leq \frac{1}{144\sqrt{2}eb_{0}^{\frac{3}{2}}\sqrt{\delta}}.
\end{eqnarray*}
These results can be summarized in the following criteria.\\
\\
{\bf Case1}. \fbox{$f\in B_{\sigma}$} Let $\sigma>0$ and $f\in B_{\sigma}$. Under the conditions of Theorem2.7, for any fixed $\delta,\ 0<\delta\leq b_{0}$, the optimal choice of $\beta$ in the interval $(0,\beta_{0}]$ where $\beta_{0}=\frac{1}{144\sqrt{2}eb_{0}^{\frac{3}{2}}\sqrt{\delta}}$ is the number minimizing
$$MN(\beta):=\beta^{\frac{b_{0}}{2\delta}}e^{\frac{\sigma^{2}}{8\beta}}.$$
{\bf Numerical Results}:\\
\\
Note that in this case the optimal choice of $\beta$ is independent of the dimension. Therefore in our numerical examples we will totally ignore the influence of n.
\begin{figure}[t]
\centering
\includegraphics[scale=1.0]{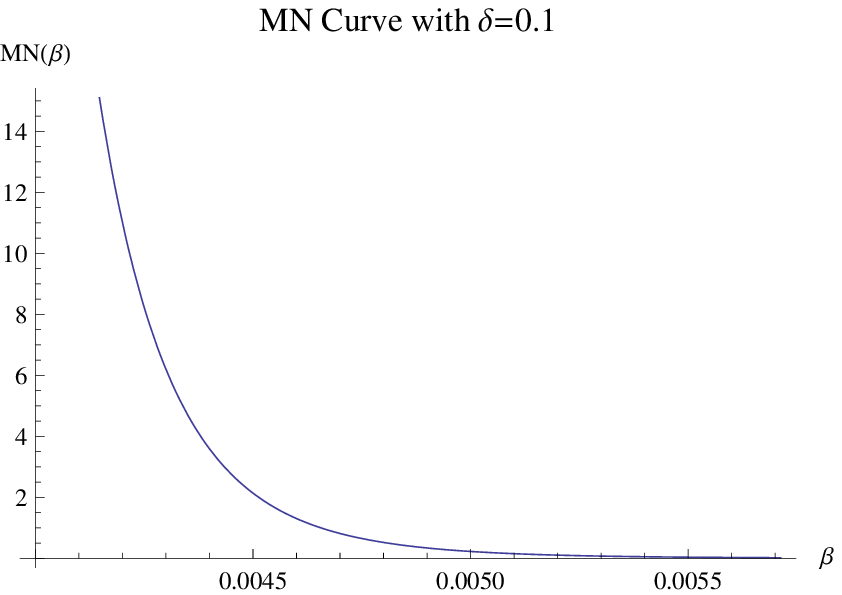}
\caption{$f\in B_{\sigma}, \sigma =1, b_{0}=1$}

\includegraphics[scale=1.0]{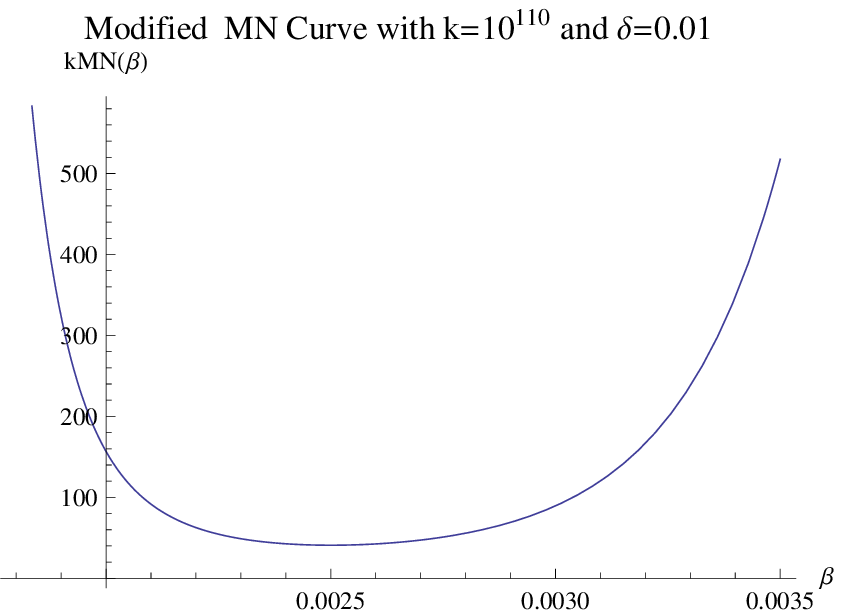}
\caption{$f\in B_{\sigma}, \sigma =1, b_{0}=1$}

\includegraphics[scale=1.0]{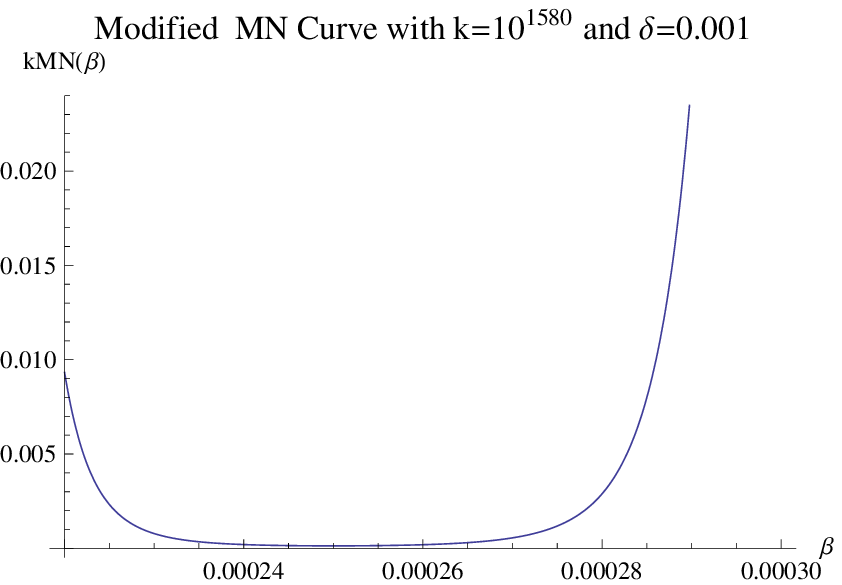}
\caption{$f\in B_{\sigma}, \sigma =1, b_{0}=1$}

\end{figure}

\clearpage

\begin{figure}[t]
\centering
\includegraphics[scale=1.0]{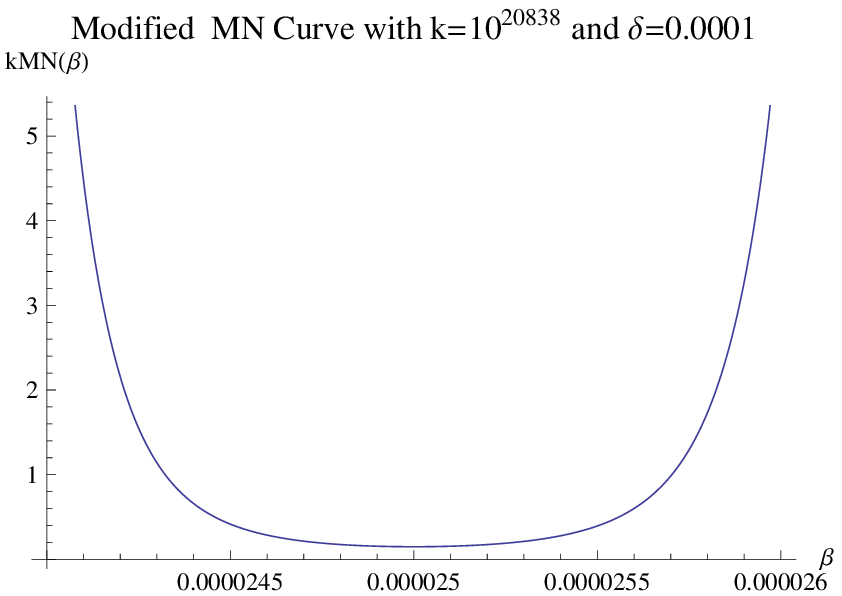}
\caption{$f\in B_{\sigma}, \sigma =1, b_{0}=1$}

\includegraphics[scale=1.0]{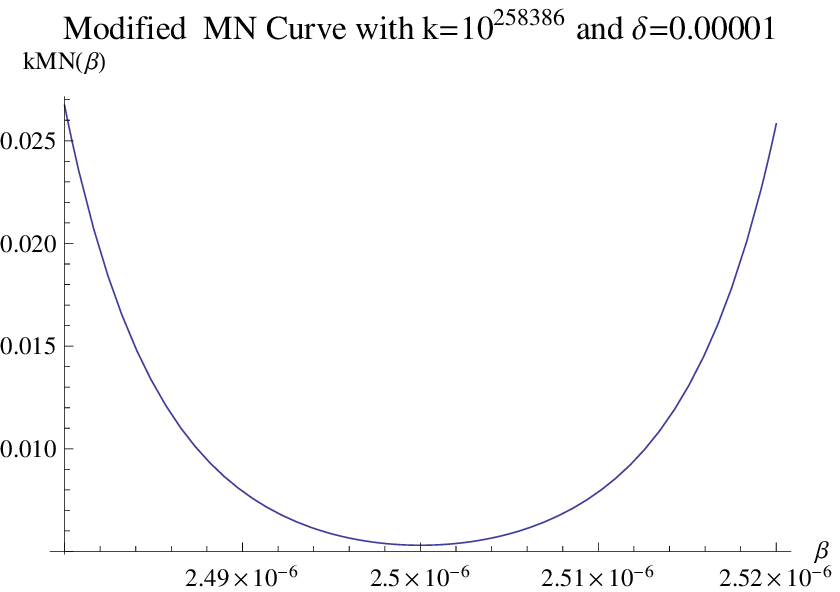}
\caption{$f\in B_{\sigma}, \sigma =1, b_{0}=1$}
\end{figure}

Now we begin the second case.\\
\\
{\bf Case2}. \fbox{$f\in G_{\sigma}$} Let $\sigma>0$ and $f\in G_{\sigma}$. Under the conditions of Theorem2.7, for any fixed $\delta,\ 0<\delta\leq b_{0}$, the optimal choice of $\beta$ in the interval $(0,\beta_{0}]$ where $\beta_{0}=\frac{1}{144\sqrt{2}eb_{0}^{\frac{3}{2}}\sqrt{\delta}}$ is the number minimizing
$$MN(\beta):=\beta^{\frac{b_{0}}{2\delta}}\left\{\sup_{\xi\in R^{n}}e^{\frac{|\xi|^{2}}{4\beta}-\frac{|\xi|^{3}}{\sigma}}\right\}^{1/2}.$$
{\bf Numerical Results}:\\
\\
Here, also, the optimal choice of $\beta$ is independent of n.
\begin{figure}[t]
\centering
\includegraphics[scale=1.1]{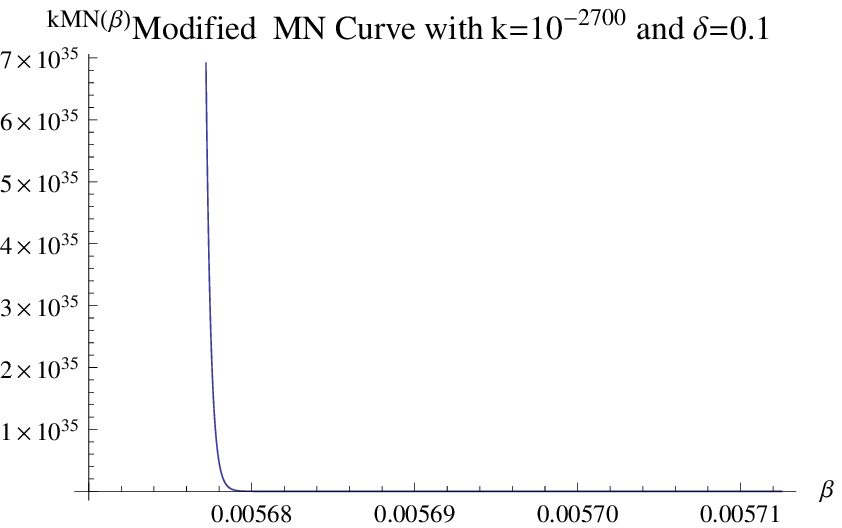}
\caption{$f\in G_{\sigma}, \sigma =1, b_{0}=1$}

\includegraphics[scale=1.0]{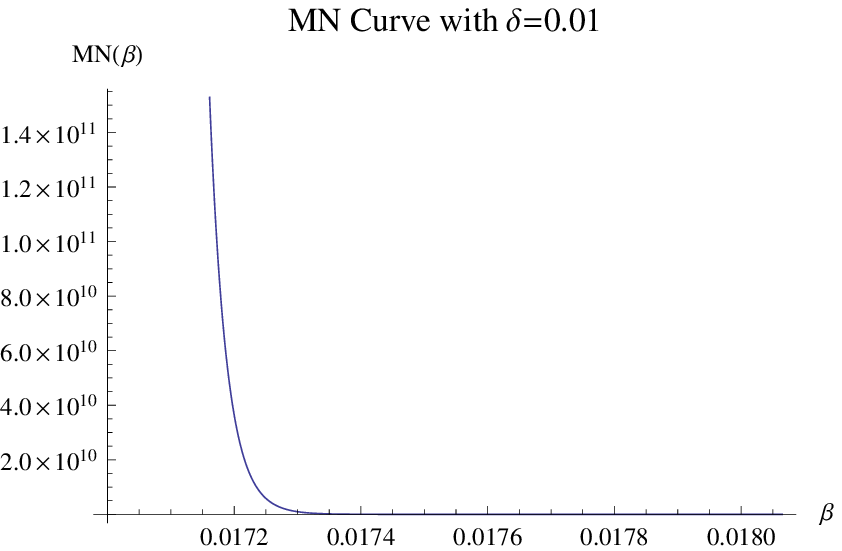}
\caption{$f\in G_{\sigma}, \sigma =1, b_{0}=1$}

\includegraphics[scale=1.1]{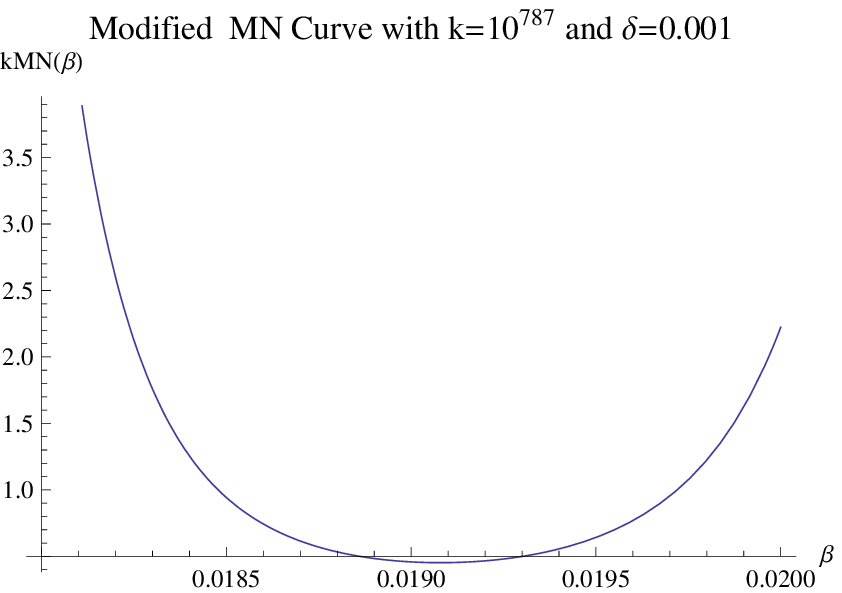}
\caption{$f\in G_{\sigma}, \sigma =1, b_{0}=1$}

\end{figure}
\clearpage

\begin{figure}[t]
\centering
\includegraphics[scale=1.0]{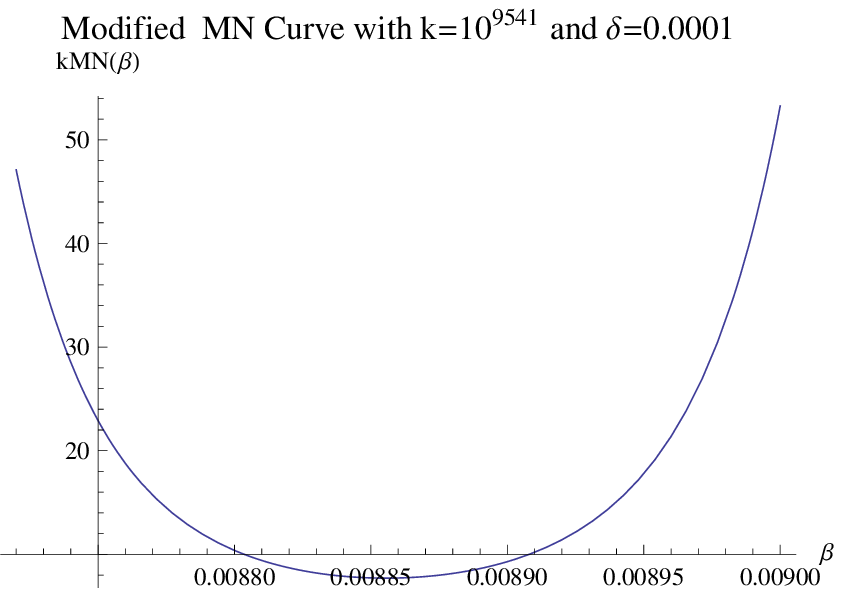}
\caption{$f\in G_{\sigma}, \sigma =1, b_{0}=1$}

\includegraphics[scale=1.0]{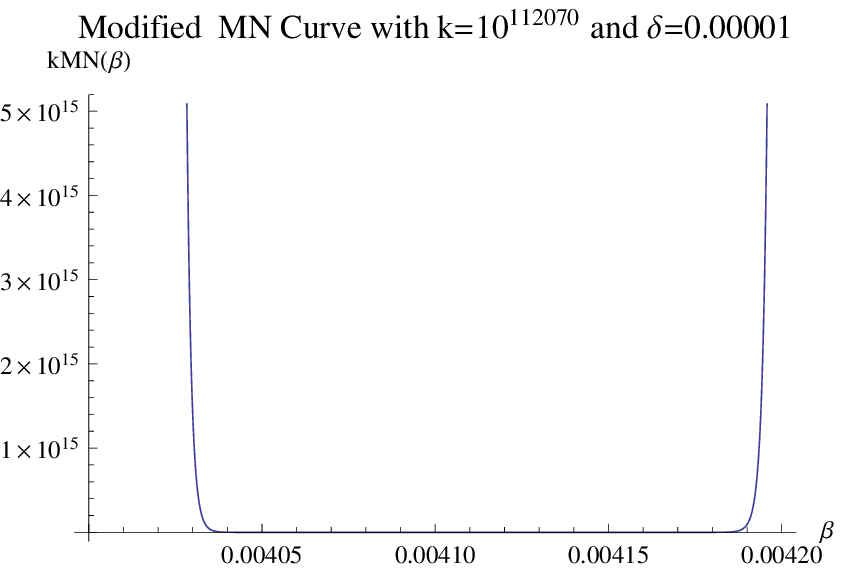}
\caption{$f\in G_{\sigma}, \sigma =1, b_{0}=1$}
\end{figure}

As in section3, although we didn't present the entire MN curves in virtue of the restriction of Mathematica, the crucial parts were presented.\\
\\
{\bf Remark}:(a)In our criteria there is an upper bound $\frac{1}{144\sqrt{2}eb_{0}^{\frac{3}{2}}\sqrt{\delta}}$ for $\beta$. This is a drawback of our theory and deserves future research. (b) Theorem2.7 requires that interpolation happens in an n-simplex with centers(interpolation points) evenly spaced points. It means that the data points are not purely scattered. This is also a drawback. However the shape of the simplex, and hence the distribution of the centers, is very flexible, making this drawback harmless. (c)In section4 the criteria do not depend on the dimension. This is the main advantage over the criteria of section3.

\end{document}